\documentclass[11pt,reqno]{amsart}

\usepackage{amssymb}
\usepackage{amsmath}
\usepackage{amsfonts}
\usepackage{amscd}
\usepackage{mathrsfs}
\usepackage{color}

\newtheorem{theorem}{Theorem}[section]

\newtheorem{lemma}[theorem]{Lemma}

\newtheorem{problem}[theorem]{Problem}

\numberwithin{equation}{section}

\begin{document}

\title[]{\small On meromorphic solutions of functional equations of  Fermat type }

\author[]{Pei-chu Hu\quad and \quad Qiong Wang*}

\address{$^{\dagger}$ Department of Mathematics, Shandong University
\vskip 2pt \hspace{1.4mm} Jinan, Shandong, 250100, P.R. China
\vskip 2pt \hspace{1.5mm} Email: {\sf pchu@sdu.edu.cn, qiongwangsdu@126.com}}

\thanks{{\sf Qiong Wang* is the corresponding author.}}
\thanks{This work of both authors was partially supported by NSFC of China (No. 11461070, 11271227), PCSIRT (No. IRT1264), and the Fundamental Research Funds of Shandong University (No. 2017JC019)}

\begin{abstract}
 Take complex numbers $a_j,b_j$, $(j=0,1,2)$ such that $c\neq0$ and
\begin{equation*}
{\rm rank}
 \left(                 
  \begin{array}{ccc}   
     a_{0} &  a_{1} &  a_{2}\\  
     b_{0} &  b_{1} &  b_{2}\\  
  \end{array}
\right)=2.                
\end{equation*}
We show that if the following functional equation of  Fermat type
\begin{equation*}
\begin{aligned}
 \left\{a_{0}f(z)+a_{1}f(z+c)+a_{2}f'(z)\right\}^3+\left\{b_{0}f(z)+b_{1}f(z+c)+b_{2}f'(z)\right\}^3=e^{\alpha z+\beta}
\end{aligned}
\end{equation*}
has meromorphic solutions of finite order, then it has only entire solutions of the form $f(z)=Ae^{\frac{\alpha z+\beta}{3}}+Ce^{Dz},$
 which generalizes the results in  \cite{19} and \cite{14}.
\end{abstract}

\maketitle

\section{Introduction}\label{Int} 
Take complex numbers $\alpha,\beta,c$ with $c\not=0$.
We will characterize all global meromorphic solutions of the following functional equation of  Fermat type
\begin{equation}
\begin{aligned}\label{eq-1.3}
  \left\{a_{0}f(z)+a_{1}f(z+c)+a_{2}f'(z)\right\}^3+\left\{b_{0}f(z)+b_{1}f(z+c)+b_{2}f'(z)\right\}^3=e^{\alpha z+\beta}
\end{aligned}
\end{equation}
under the assumption:
\begin{description}
\item[(A)] Take six complex numbers $a_i,b_i$ for $i=0,1,2$ such that
\begin{equation*}
{\rm rank}
 \left(                 
  \begin{array}{ccc}   
     a_{0} &  a_{1} &  a_{2}\\  
     b_{0} &  b_{1} &  b_{2}\\  
  \end{array}
\right)=2.                
\end{equation*}
 \end{description}

This kind of problems goes back to classical functional equation of Fermat type
\begin{equation}\label{eq-1.1}
f^{n}+g^{n}=1
\end{equation}
for a positive integer $n$. When $n>3$, Gross \cite{11} proved that the equation $(\ref{eq-1.1})$ has no  nonconstant meromorphic solutions on the complex plane $\mathbb{C}$. If $n=2$, Gross \cite{11} showed that all meromorphic solutions of $(\ref{eq-1.1})$ on $\mathbb{C}$ are of the form
\begin{equation*}
  f=\frac{2\omega}{1+\omega^{2}},\ \   g=\frac{1-\omega^{2}}{1+\omega^{2}},
\end{equation*}
 where $\omega$ is a nonconstant meromorphic function on $\mathbb{C}$. If $n=3$, the form of meromorphic solutions of $(\ref{eq-1.1})$ was conjectured by Gross \cite{11}, and was completely characterized by Baker \cite{15} (see Section~\ref{Sec2}). Meanwhile, we refer the reader
 to  \cite{12}, \cite{13} and \cite{18}.

Take two positive integers $m$ and $n$. When $m\geq3$ and $n\geq3$, Montel \cite{18} proved that the following functional equation
\begin{equation}\label{m-n-equ}
f^{n}+g^{m}=1
\end{equation}
has no transcendental entire solutions. Further, Yang \cite{7} showed that the equation (\ref{m-n-equ}) has no nonconstant entire solutions if $\lambda=1-\frac{1}{m}-\frac{1}{n}>0$. For more detail, we refer the reader to the work of Hu, Li and Yang \cite{26}.
For some works related to partial differential equations of Fermat type, see \cite{10}, \cite{20}, \cite{2}-\cite{5}.

In 1985, Hayman \cite{28} proved that the following functional equation
\begin{equation}\label{n3equ}
f^{n}+g^{n}+h^{n}=1
\end{equation}
has no nonconstant meromorphic (resp., entire) solutions if $n\geq9$ (resp., $ n\geq6$).
For the case $n=5$ and $n=6$,
Gundersen \cite{40}, \cite{41} proved the existence of transcendental meromorphic solutions of the equation (\ref{n3equ}).
When $n=8$, Ishizaki \cite{27}  showed that if $f, g, h$ are
non-constant meromorphic functions satisfying (\ref{n3equ}), then there must
exist a small function $a(z)$ with respect to $f, g,$ $h$ such that
$$W( f^{8}, g^{8}, h^{8}) = a(z)(f(z)g(z)h(z))^{6},$$
where the left side is a Wronskian.
Ng and Yeung \cite{29} proved that (\ref{n3equ}) has no non-trivial meromorphic (resp. entire) solutions when $n\geq8 $ (resp. $n\geq6$).
Gundersen \cite{35} ask that do (\ref{n3equ}) have non-constant meromorphic
solutions for $n=7$ or $n=8$?
Futhermore, Ng and Yeung \cite{29}, Yang \cite{42} studied the existence of meromorphic solutions of the equation
$$f^{n}+g^{m}+h^{p}=1,$$
where $n, m, p$ are positive integers.

More general, if $k(\geq3)$ positive integers $n_1,...,n_k$ satisfy
\begin{equation*}
    \lambda=1-\sum_{j=1}^k\frac{k-1}{n_j}>0,
\end{equation*}
Hu, Li and Yang \cite{26} proved that the following functional equation
\begin{equation}\label{k-thequ}
  f_{1}^{n_{1}}+f_{2}^{n_{2}}+\cdots +f_{k}^{n_{k}}=1
\end{equation}
has no nonconstant entire solutions.

Based on the observation above, these functional equations maybe have global nonconstant meromorphic solutions when the powers are lower, which further are characterized by some researchers. Hence, a natural question follows:
\begin{problem}
When the powers $n_j$ in the equation $(\ref{k-thequ})$ are lower, can we characterize all meromorphic solutions $f$ of the equation $(\ref{k-thequ})$ if $f_j$ are replaced by differentials or differences $($even their mixture$)$ of a fixed function $f$?
\end{problem}

For example, we may try to characterize all meromorphic solutions $f$ of the following difference equation
\begin{equation}\label{nm-dif-equ}
f^{n}(z)+f^{m}(z+c)=1
\end{equation}
for a fixed non-zero constant $c$ (cf. \cite{18}, \cite{8}).
In fact, Shimomura \cite{23} considered the case $n>m=1$ and proved that there always exists an entire solutions of infinite order. Later, Liu \cite{16} investigated that difference equation (\ref{nm-dif-equ}) has no transcendental entire solutions with finite order.
 Liu et al. \cite{17} illustrated that the solutions of (\ref{nm-dif-equ}) are periodic functions with period $2c$ for $n=m=1,$ and found that (\ref{nm-dif-equ}) has transcendental entire solutions with finite order for $n=m=2$.
 
 In 2004, Yang and Li \cite{9} proved that the following differential equation
\begin{equation*}
f^{2}(z)+(f'(z))^{2}=1
\end{equation*}
has only transcendental entire solutions of the form
\begin{equation*}
 f(z)=\frac{1}{2}\left(Pe^{\alpha z}+\frac{1}{P}e^{-\alpha z}\right),
\end{equation*}
where $P, \alpha$ are nonzero constants. However, the following differential equation
\begin{equation*}
f^{3}(z)+(f'(z))^{3}=1
\end{equation*}
has no nonconstant global meromorphic solutions (cf. \cite{14}). For more results related to shifts and differences of entire and
meromorphic functions, see \cite{25}, \cite{32}, \cite{31}, \cite{43}, \cite{16}, \cite{17}, \cite{23}.

Recently, Han and L\"{u} \cite{19} studied that the following difference equation
\begin{equation*}
    f^{3}(z)+f^{3}(z+c)=e^{\alpha z+\beta}
\end{equation*}
has only entire solutions of the form $f(z)=Ae^{\frac{\alpha z+\beta}{3}}$
with a constant $A$ satisfying $A^{3}(1+e^{\alpha c})=1$.

In this paper, we will discuss the existence of solutions to functional equations of Fermat type, and shall obtain the following main result.
\begin{theorem}\label{mai-1.4}
Take complex numbers $\alpha,\beta,c,a_i,b_i,i=0,1,2$ with $c\not=0$, and assume $({\rm A})$.
If the equation $(\ref{eq-1.3})$ has meromorphic solutions of finite order, then
it has only entire solutions of the following form
\begin{equation*}
  f(z)=Ae^{\frac{\alpha z+\beta}{3}}+Ce^{Dz},
\end{equation*}
where $A,C,D$ are constants. Moreover, if we define constants $c_{0},c_{1}$ by $c_{0}^{3}+c_{1}^{3}=1$, then $A,D$ are completely determined by $a_{i},b_{i},\alpha, c_{0},c_{1}$ as follows:
\begin{enumerate}
  \item  $A=\frac{b_{2}c_{0}-a_{2}c_{1}}{a_{0}b_{2}-a_{2}b_{0}}$, $C=0$ if $a_{0}b_{1}-a_{1}b_{0}=0$, $b_{1}a_{2}-a_{1}b_{2}=0$,
  \item  $A=\frac{b_{1}c_{0}-a_{1}c_{1}}{a_{0}b_{1}-a_{1}b_{0}}$, $C=0$ if $a_{0}b_{1}-a_{1}b_{0}\neq0$, $b_{1}a_{2}-a_{1}b_{2}=0$,
  \item  $A=\frac{3(b_{1}c_{0}-a_{1}c_{1})}{(b_1a_2-a_1b_2)(\alpha-3D)}$, $C\in\mathbb{C}$, $D=\frac{a_{1}b_{0}-a_{0}b_{1}}{b_{1}a_{2}-a_{1}b_{2}}$ if $b_{1}a_{2}-a_{1}b_{2}\neq0$.
\end{enumerate}
\end{theorem}

In particular, if we take $a_{0}=1$,~$a_{1}=a_{2}=0$, $b_{1}=1$,~$b_{0}=b_{2}=0$ in Theorem \ref{mai-1.4}, this is just the results in \cite{19} and \cite{14}.

\section{Preliminary}\label{Sec2}

We assume that the reader is familiar with Nevanlinna theory (cf. \cite{45}, \cite{21}) of meromorphic functions $f$ in $\mathbb{C}$, such as the first main theorem of $f$, the second main theorem of $f$, the characteristic function $T(r,f)$, the proximity function $m(r,f)$, the counting functions $N(r,f)$, $\bar{N}(r,f)$, and $S(r,f)$, where as usual $S(r,f)$ denotes any quantity satisfying $S(r,f)=o(T(r,f))$ as $r\rightarrow\infty$ outside a possible exceptional set of finite logarithmic measure. Further, recall that the order of $f$ is defined by
\begin{equation*}
    \rho(f)=\lim \limits_{r\rightarrow\infty}\sup\frac{\log T(r,f)}{\log r}.
\end{equation*}

When $n=3$ in the equation $(\ref{eq-1.1})$, Gross  \cite{11} showed that the equation $(\ref{eq-1.1})$ has a pair of meromorphic solutions
\begin{equation}\label{gross-fun}
  f(z)=\frac{1}{2 \wp(z)}\left\{1+\frac{\wp'(z)}{\sqrt{3}}\right\}, ~~g(z)=\frac{1}{2\wp(z)}\left\{1-\frac{\wp'(z)}{\sqrt{3}}\right\},
\end{equation}
where $\wp(z)$ is the Weierstrass elliptic function with periods $\omega_{1}$ and $\omega_{2}$, defined by
\begin{equation*}
  \wp (z;\omega_{1},\omega_{2}):=\frac{1}{z^{2}}+\sum \limits _{\mu,\nu\in\mathbb{Z}, \mu^{2}+\nu^{2}\neq0 }\left\{
  \frac{1}{(z+\mu\omega_{1}+\nu\omega_{2})^{2}}-\frac{1}{(\mu\omega_{1}+\nu\omega_{2})^{2}}\right\},
\end{equation*}
which is an even function and satisfies, after appropriately choosing $\omega_{1}$ and $\omega_{2}$,
\begin{equation}\label{eq-1.2}
  (\wp')^{2}=4\wp^{3}-1.
\end{equation}
Furthermore, Baker \cite{15} provided the following important result.

\begin{lemma}[\cite{15}]\label{th-1.1}
Any nonconstant meromorphic functions $F(z),G(z)$ in
the complex plane $\mathbb{C}$ satisfying
\begin{equation}\label{EqF-G}
  F^{3}(z)+G^{3}(z)=1
\end{equation}
have the form $F(z)=f(h(z)), G(z)=\eta g(h(z))=\eta f(-h(z))$, where $f$ and $g$ are elliptic functions defined by $(\ref{gross-fun})$, $h(z)$ is an entire function in $\mathbb{C}$ and $\eta$ is a cube-root of the unity.
\end{lemma}

The following lemma is referred to Bergweiler \cite{24} and Edrei and Fuchs \cite{1}.

\begin{lemma}[\cite{24,1}]\label{le1}
Let $f$ be a meromorphic function and $h$ be an entire function in $\mathbb{C}$. When $0<\rho(f),$ $\rho(h)<\infty$, then $\rho(f\circ h)=\infty$. When $\rho(f\circ h)<\infty$
and $h$ is transcendental, then $\rho(f)=0$.
\end{lemma}

\begin{lemma}[\cite{26}]\label{le2}
If $f$ is a nonconstant meromorphic function, then
\begin{equation*}
  T(r,A_{f})=pT(r,f)+O\left(\sum \limits^{p}_{j=0}T(r,a_{j})\right),
\end{equation*}
where $a_{j}$ are meromorphic functions with $a_{p}\not\equiv0$ and
\begin{equation*}
 A_{f}=A(z,f(z))=\sum \limits^{p}_{j=0}a_{j}(z)f^{j}(z).
\end{equation*}
\end{lemma}

\begin{lemma}[\cite{25}]\label{le3}
Let $f(z)$ be a meromorphic function such that the order $\rho(f)<+\infty$,
and let $c$ be a nonzero complex number. Then for each $\varepsilon >0$, we have
\begin{equation*}
  T(r,f_c) = T(r,f)+O(r^{\rho(f) -1+\varepsilon})+O(\log r),
\end{equation*}
where $f_c(z)=f(z+c)$.
\end{lemma}

\section{Proof of Theorem~\ref{mai-1.4}}
First of all, we rewrite (\ref{eq-1.3}) into the form (\ref{EqF-G}), where $F$ and $G$ are defined by
\begin{equation}\label{eq-new1-1}
\begin{aligned}
  &F(z)e^{\frac{\alpha z+\beta}{3}}=a_{0}f(z)+a_{1}f(z+c)+a_{2}f'(z),\\
  &G(z)e^{\frac{\alpha z+\beta}{3}}=b_{0}f(z)+b_{1}f(z+c)+b_{2}f'(z).
  \end{aligned}
\end{equation}
Then we claim that $F(z),G(z)$ are constants.

We assume, to the contrary, that $F(z),G(z)$ are not constants. By Lemma \ref{th-1.1}, we have
\begin{equation}\label{eq-new1}
F(z)=\frac{1+\frac{\wp'(h(z))}{\sqrt{3}}}{2\wp(h(z))},\ \ G(z)=\frac{1-\frac{\wp'(h(z))}{\sqrt{3}}}{2\wp(h(z))}\eta .
\end{equation}
Based on ideas in \cite{19} and \cite{14}, we confirm a fact as follows:

\textbf{Claim 1:} $h(z)$ must be a nonconstant polynomial.

Solving $\wp'(h(z))$ from the first expression of (\ref{eq-new1}), we obtain
\begin{equation*}
\wp'(h(z))=\sqrt{3} (2F(z)\wp(h(z))-1).
\end{equation*}
Substituting into (\ref{eq-1.2}), it follows that
\begin{equation}\label{eq-2.16}
\wp^{3}(h(z))=  3F^{2}(z)\wp^{2}(h(z))-3F(z)\wp(h(z))+1,
\end{equation}
which immediately implies the following inequality of Nevanlinna's characteristic functions
\begin{equation*}
  3T(r,\wp(h))\leq 2T(r,F)+2T(r,\wp(h))+O(1),
\end{equation*}
that is, $T(r,\wp(h))\leq 2T(r,F)+O(1).$
Hence we have $\rho(\wp(h))\leq \rho(F)$. The first expression of (\ref{eq-new1-1}) means that $\rho(F)<\infty$ since $\rho(f)<\infty$ by the assumption. Therefore, we obtain $\rho(\wp(h))<\infty$.

Bank and Langley \cite[equation (2.7)]{22} gave Nevanlinna's characteristic function of $\wp$ as follows:
\begin{equation}\label{eq-2.17}
  T(r,\wp)=\frac{\pi}{A}r^{2}(1+o(1)) ,
\end{equation}
where $A$ is the area of a parallelogram $P_{a}$
with four vertices $0,\omega_{1},\omega_{2},\omega_{1}+\omega_{2}$, which further yields $\rho(\wp)=2$. Then  $h$ must be a polynomial based on Lemma \ref{le1}. Further, by using(\ref{eq-new1}), we know that $h(z)$ is not a constant.

Next we distinguish three cases to prove Theorem~\ref{mai-1.4}.

\textbf{Case 1:} $a_{0}b_{1}-a_{1}b_{0}=0$, $b_{1}a_{2}-a_{1}b_{2}=0$.

Then we have $a_{0}b_{2}-b_{0}a_{2}\neq0$ by the rank assumption. Solving the equations (\ref{eq-new1-1}) and noting (\ref{eq-new1}), it follows that
\begin{equation}\label{eq-new2}
  f(z)=\frac{b_{2}F(z)-a_{2}G(z)}{a_{0}b_{2}-b_{0}a_{2}}e^{\frac{\alpha z+\beta}{3}}=\frac{b_{2}-a_{2}\eta+\frac{b_{2}+a_{2}\eta}{\sqrt{3}}\wp'(h(z))}{2(a_{0}b_{2}-b_{0}a_{2})\wp(h(z))}e^{\frac{\alpha z+\beta}{3}}.
\end{equation}
By differentiating (\ref{eq-new2}) and noting that $(\wp')^{2}=4\wp^{3}-1$, $\wp''=6\wp^{2}$, we have
\begin{equation}\label{eq-new3}
\begin{aligned}
f'(z)=&\frac{\frac{2(b_{2}+a_{2}\eta)}{\sqrt{3}}\wp^{3}(h(z))h'(z)-(b_{2}-a_{2}\eta)\wp'(h(z))h'(z)+\frac{b_{2}+a_{2}\eta}{\sqrt{3}}h'(z)}{2(a_{0}b_{2}-b_{0}a_{2})\wp^{2}(h(z))}e^{\frac{\alpha z+\beta}{3}}\\
&+\frac{\alpha\left(b_{2}-a_{2}\eta+\frac{b_{2}+a_{2}\eta}{\sqrt{3}}\wp'(h(z))\right)}{6(a_{0}b_{2}-b_{0}a_{2})\wp(h(z))}e^{\frac{\alpha z+\beta}{3}}.
\end{aligned}
\end{equation}
Substituting (\ref{eq-new2}) and (\ref{eq-new3}) into the first equation of (\ref{eq-new1-1}), we have
\begin{equation}\label{eq-new4}
\begin{aligned}
&\Bigg((a_{0}b_{2}-b_{0}a_{2})\Big(1+\frac{\wp'(h(z))}{\sqrt{3}}\Big)\wp(h(z))-\frac{a_{0}+a_{2}\alpha}{3}\Big(b_{2}-a_{2}\eta
+\frac{b_{2}+a_{2}\eta}{\sqrt{3}}\wp'(h(z))\Big)\\
&\wp(h(z))-a_{2}\Big(\frac{2(b_{2}+a_{2}\eta)}{\sqrt{3}}\wp^{3}(h(z))h'(z)
-(b_{2}-a_{2}\eta)\wp'(h(z))h'(z)+\frac{b_{2}+a_{2}\eta}{\sqrt{3}}h'(z)\Big)\Bigg)\\
&\wp(h(z+c))=e^{\frac{\alpha c}{3}}a_{1}\left(b_{2}-a_{2}\eta+\frac{b_{2}+a_{2}\eta}{\sqrt{3}}\wp'(h(z+c))\right)\wp^{2}(h(z)).
 \end{aligned}
\end{equation}

Let $\{z_{j}\}^{\infty}_{j=1}$ be the zeroes of $\wp$  satisfying $|z_{j}|\rightarrow \infty$ as $j\rightarrow \infty$. Then for each $j$,  there exist $\deg(h)$ complex numbers $a_{j,k}$ such that $h(a_{j,k})=z_j$, for $\ k=1,...,\deg(h)$.
Moreover, the equation (\ref{eq-1.2}) implies $(\wp')^{2}(h(a_{j,k}))=(\wp')^{2}(z_{j})=-1$, since $\wp(z_{j})=0$.

Now, we confirm the second fact as follows:

\textbf{Claim 2:} $\wp(h(a_{j,k}+c))=0$ only holds for at most finitely many $a_{j,k}$'s.

We assume, to the contrary, that there exists an infinite subsequence of  $\{a_{j,k}\}$, without loss of generality we may take $\{a_{j,k}\}$  itself,  such that $\wp(h(a_{j,k}+c))=0$. Thus we have
$(\wp')^{2}(h(a_{j,k}+c))=-1$.
 Differentiating (\ref{eq-new4}) and then setting $z=a_{j,k}$, we derive a relation
 \begin{equation}\label{new-1-2}
   a_2 h'(a_{j,k})h'(a_{j,k}+c)\left((a_{2}\eta-b_2)\wp'(h(a_{j,k}))+\frac{b_{2}+a_{2}\eta}{\sqrt{3}}\right)=0.
 \end{equation}
Note that $\wp'(h(a_{j,k}))=\pm i$, where $i$ is the imaginary unit. If $a_{2}\neq0$, the equation (\ref{new-1-2}) becomes
\begin{equation*}
  h'(a_{j,k})h'(a_{j,k}+c)\left(\pm(a_{2}\eta-b_2)i+\frac{b_{2}+a_{2}\eta}{\sqrt{3}}\right)=0.
\end{equation*}
Since $h$ is a nonconstant polynomial and $a_{j,k}$ is an infinite sequence, then there are infinitely many $a_{j,k}$ such that  $h'(a_{j,k})h'(a_{j,k}+c)\neq 0$. It follows that
\begin{equation*}
    (a_{2}\eta-b_2)i+\frac{b_{2}+a_{2}\eta}{\sqrt{3}}=\frac{1-\sqrt{3}i}{\sqrt{3}}(b_2+e^{\frac{2\pi i}{3}}\eta a_2)=0
\end{equation*}
or
\begin{equation*}
 (b_{2}-a_{2}\eta)i+\frac{b_{2}+a_{2}\eta}{\sqrt{3}}=\frac{1+\sqrt{3}i}{\sqrt{3}}(b_2+e^{-\frac{2\pi i}{3}}\eta a_2)=0,
\end{equation*}
and hence $b_2+e^{\frac{\pm2\pi i}{3}}\eta a_2=0.$
We may let $b_2=-k a_2,$
where $k=e^{\frac{\pm2\pi i}{3}}\eta$.
Note that $a_{2}\neq0$, then we get $ b_1=-k a_1$
because the assumption $b_{1}a_{2}-a_{1}b_{2}=0$.
Combining with the assumption $a_{0}b_{1}-a_{1}b_{0}=0$, we obtain $ka_{0}a_{1}+a_{1}b_{0}=0.$

If $a_{1}\neq0$, $ka_{0}a_{1}+a_{1}b_{0}=0$ implies $ b_0=-k a_0.$
It follows that $b_i=-ka_i$ for $i=0,1,2,$
which contradicts to the rank assumption.

When $a_{1}=0$, by using (\ref{eq-new1}), we can rewrite the first equation of (\ref{eq-new1-1}) as follows
\begin{equation}\label{eq-new-102}
  f'(z)=-\frac{a_{0}}{a_{2}}f(z)+\frac{1+\frac{\wp'(h(z))}{\sqrt{3}}}{2a_2\wp(h(z))}e^{\frac{\alpha z+\beta}{3}}.
\end{equation}
By using the theory of first-order linear differential equations, the solution of (\ref{eq-new-102}) is
\begin{equation}\label{eq-new103}
  f(z)=e^{-\frac{a_{0}}{a_{2}}z}\left(\int\frac{1+\frac{\wp'(h(z))}{\sqrt{3}}}{2a_2\wp(h(z))}
  e^{\frac{3a_{0}+\alpha a_{2}}{3a_{2}}z+\frac{\beta}{3}}dz+C\right),
\end{equation}
where $C$ is a constant. Let $t_j\ (j\geq 1)$ be the poles of $\wp(z)$. Then we have $ \wp(z)=g_j(z)(z-t_{j})^{-2},$
where $g_j$ is a holomorphic function in a neighborhood of $t_j$ with $g_j(t_j)\not=0$, and hence
\begin{equation*}
\frac{1+\frac{\wp'(z)}{\sqrt{3}}}{2a_{2}\wp(z)}=-\frac{1}{\sqrt{3}a_{2}(z-t_{j})}
+\frac{(z-t_{j})^{2}}{2a_{2}g_j(z)}+\frac{g_j'(z)}{2\sqrt{3}a_{2}g_j(z)}
=-\frac{1}{\sqrt{3}a_{2}(z-t_{j})}+O(1)
\end{equation*}
as $z\rightarrow t_j$.
Setting $z=h(z)$, we derive a relation
\begin{equation*}
\begin{aligned}
\frac{1+\frac{\wp'(h(z))}{\sqrt{3}}}{2a_{2}\wp(h(z))}=-\frac{1}{\sqrt{3}a_{2}(h(z)-t_{j})}+O(1)
\end{aligned}
\end{equation*}
as $h(z)\rightarrow t_j$. Note that the equation $h'(z)=0$ only has finitely many solutions, but $\bigcup \limits^{\infty}_{j=1}h^{-1}(t_{j})$ is an infinite set. Thus there exist an integer $j$ and a point $z'\in h^{-1}(t_{j})$ such that $h'(z')\neq 0$, that is, $z'$ is a simple zero of $h(z)-t_j$. Therefore $f(z)$ has a logarithmic singular point $z'$, which contradicts the assumption  that $f(z)$ is a meromorphic function.

Next we consider the case $a_2=0$. Now we claim that $a_{1}=0$. Otherwise, if $a_{1}\neq0$, the assumption $b_{1}a_{2}-a_{1}b_{2}=0$ implies $b_{2}=0$, which is contradicted to the rank assumption $b_{2}a_{0}-a_{2}b_{0}\neq0$.
Since $a_1=a_2=0$, we find $b_{2}a_{0}\neq0$
from $b_{2}a_{0}-a_{2}b_{0}\neq0$ and $a_{0}b_{1}=0$ because the assumption $a_{0}b_{1}-a_{1}b_{0}=0$. Thus it follows that $a_{0}\neq0$, $b_{1}=0$, $b_{2}\neq0$.
By using (\ref{eq-new1}), we can rewrite the equations (\ref{eq-new1-1}) as follows
\begin{equation}\label{eq-new5}
f(z)=\frac{1+\frac{\wp'(h(z))}{\sqrt{3}}}{2a_0\wp(h(z))}e^{\frac{\alpha z+\beta}{3}},\
b_{0}f(z)+b_{2}f'(z)=\frac{1-\frac{\wp'(h(z))}{\sqrt{3}}}{2\wp(h(z))}\eta e^{\frac{\alpha z+\beta}{3}}.
\end{equation}
Differentiating the first relation of (\ref{eq-new5}) and noticing that $(\wp')^{2}=4\wp^{3}-1$, $\wp''=6\wp^{2}$, it follows that
\begin{equation}\label{eq-new7}
f'(z)=\frac{\frac{2}{\sqrt{3}}\wp^{3}(h(z))h'(z)-\wp'(h(z))h'(z)+\frac{h'(z)}{\sqrt{3}}}{2a_{0}\wp^{2}(h(z))}e^{\frac{\alpha z+\beta}{3}}
+\frac{\alpha\left(1+\frac{\wp'(h(z))}{\sqrt{3}}\right)}{6a_{0}\wp(h(z))}e^{\frac{\alpha z+\beta}{3}}.
\end{equation}
Substituting $f$ and $f'$ into the second equation of (\ref{eq-new5}), we have
\begin{equation*}
\begin{aligned}
&\Bigg(\frac{2b_{2}}{\sqrt{3}}\wp^{3}(h(z))h'(z)-\Big(a_{0}\eta-b_{0}-\frac{b_{2}\alpha}{3}\Big)\wp(h(z))+ \frac{ b_{2}}{\sqrt{3}}h'(z)\Bigg)^{2}\\
 =&(4\wp^3(h(z))-1)\Bigg(b_{2}h'(z)-\frac{a_{0}\eta+b_{0}+\frac{b_{2}\alpha}{3}}{\sqrt{3}}\wp(h(z))\Bigg)^{2}.
 \end{aligned}
\end{equation*}
Therefore, by using Lemma \ref{le2}, we find
\begin{equation*}
    6T(r,\wp(h))+O(\log r)=5T(r,\wp(h))+O(\log r),
\end{equation*}
 that is $T(r,\wp(h))=O(\log r),$
which means that $\wp(h)$ is a rational function. This is a contradiction because $\wp(h)$ is a transcendental meromorphic function. Hence Claim 2 is proved.

Based on Claim 2 and noting that $h'$ has only finitely many zeroes, there exists a positive integer $J$ such that
\begin{equation*}
    \wp(h(a_{j,k}+c))\neq0,\ h'(a_{j,k})\not=0, \ j>J,\ k=1,...,\deg(h).
\end{equation*}
Note that $a_2\not=0$, otherwise, we can obtain a contradiction according to the above arguments. Now, return to the equation (\ref{eq-new4}). We see that the coefficient of $\wp(h(z+c))$ at $a_{j,k}$ takes a nonzero value
\begin{equation*}
     a_2 h'(a_{j,k})\left((a_{2}\eta-b_2)\wp'(h(a_{j,k}))+\frac{b_{2}+a_{2}\eta}{\sqrt{3}}\right)\not=0,
\end{equation*}
when $j>J, \ k=1,...,\deg(h)$, because
\begin{equation*}
    (a_{2}\eta-b_2)\wp'(h(a_{j,k}))+\frac{b_{2}+a_{2}\eta}{\sqrt{3}}\not=0,
\end{equation*}
see the proof of Claim 2. Therefore, the equation (\ref{eq-new4}) valued at $a_{j,k}$ immediately yields
\begin{equation*}
 \wp(h(a_{j,k}+c))=\infty,    \ j>J,\ k=1,...,\deg(h),
\end{equation*}
which further yields
\begin{equation*}
 \bar{N}\left(r,\frac{1}{\wp(h)}\right)\leq \bar{N}\left(r,\wp(h_c)\right) +S(r,\wp(h)),
\end{equation*}
where $h_c(z)=h(z+c)$.

Note that the multiple zeros of $\wp(h)$ occur at zeros of its derivative $\{\wp(h)\}'=\wp'(h)h'$, that is, the zeros of $h'$ because $\wp'(h(a_{j,k}))=\pm i\not=0$. Hence we obtain an estimate
\begin{equation}\label{eq-2.21}
\begin{aligned}
  N\left(r,\frac{1}{\wp(h)}\right)&\leq \bar{N}\left(r,\frac{1}{\wp(h)}\right)+2N\left(r,\frac{1}{h'}\right) \\
  &\leq \bar{N}\left(r,\wp(h_c)\right)+S(r,\wp(h)).
  \end{aligned}
\end{equation}

Further, we claim
\begin{equation*}
    m(r,F)=S(r,\wp(h)).
\end{equation*}
In fact, the expression of $F$ in (\ref{eq-new1}) yields
\begin{equation*}
  T(r,F)\leq T(r,\wp(h))+T(r,\wp'(h))+O(1)\leq O(T(r,\wp (h))).
\end{equation*}
Based on the factorization $-G^{3}=F^{3}-1=(F-1)(F-\eta )(F-\eta^{2}),$
where $\eta\not=1$, we see that all zeros of $F-1,F-\eta $ and $F-\eta^{2} $ are of multiplicities $\geq 3$. Hence Nevanlinna's main theorems give
\begin{equation*}
\begin{aligned}
  2T(r,F)&\leq \sum \limits^{3}_{m=1}\bar{N}\left(r,\frac{1}{F-\eta^{m}}\right)+\bar{N}(r,F)+S(r,F)\\
  &\leq \frac{1}{3}\sum \limits^{3}_{m=1}N\left(r,\frac{1}{F-\eta^{m}}\right)+N(r,F)+S(r,F)\\
  & \leq T(r,F)+N(r,F)+S(r,\wp(h)),
  \end{aligned}
\end{equation*}
which immediately implies the claim $  m(r,F)=S(r,\wp(h)).$

Now we rewrite (\ref{eq-new1}) into the following form
\begin{equation*}
    \frac{1}{\wp(h)}=2F-\frac{\wp'(h)}{\sqrt{3}\wp(h)},
\end{equation*}
which means
\begin{equation*}
m\left(r,\frac{1}{\wp(h)}\right) \leq m(r,F)+ m\left(r,\frac{\wp'(h)}{\wp(h)}\right)+O(1).
\end{equation*}
Applying the lemma of logarithmic derivative, we have
\begin{equation*}
m\left(r,\frac{\wp'(h)}{\wp(h)}\right)\leq m\left(r,\frac{\wp'(h)h'}{\wp(h)}\right)    +m\left(r,\frac{1}{h'}\right) =S(r,\wp(h)),
\end{equation*}
and hence
\begin{equation}\label{eq-2.23}
m\left(r,\frac{1}{\wp(h)}\right)=S(r,\wp(h)).
\end{equation}
Combining (\ref{eq-2.21}) with (\ref{eq-2.23}), and noticing that each pole of $\wp(z)$ has multiplicity 2,  then Nevanlinna's first main theorem implies
\begin{equation*}
\begin{aligned}
  T(r,\wp(h))&=T\left(r,\frac{1}{\wp(h)}\right)+O(1)=m\left(r,\frac{1}{\wp(h)}\right)+N\left(r,\frac{1}{\wp(h)}\right)+O(1)\\
  &\leq \bar{N}(r,\wp(h_c))+S(r,\wp(h))\leq\frac{1}{2}N(r,\wp(h_c))+S(r,\wp(h))\\
  &\leq\frac{1}{2}T(r,\wp(h_c))+S(r,\wp(h))\leq\frac{1}{2}T(r,\wp(h))+S(r,\wp(h)),
  \end{aligned}
\end{equation*}
where Lemma~\ref{le3} was applied. We obtain a contradiction again.

Hence $F(z)$ and $G(z)$ are constants. We assume that $F(z)=c_{0}, G(z)=c_{1}$, where $c_{0},c_{1}$ are constants with $c_{0}^{3}+c_{1}^{3}=1$. By (\ref{eq-new2}), we have the solution
\begin{equation*}
  f(z)=\frac{b_{2}c_{0}-a_{2}c_{1}}{a_{0}b_{2}-b_{0}a_{2}}e^{\frac{\alpha z+\beta}{3}}.
\end{equation*}
This proves Case $1$ in Theorem  \ref{mai-1.4}.

\textbf{Case 2:} $a_{0}b_{1}-a_{1}b_{0}\neq0$, $b_{1}a_{2}-a_{1}b_{2}=0$.

Solving the equations (\ref{eq-new1-1}) and noting (\ref{eq-new1}), it follows that
\begin{equation}\label{eq-new8}
  f(z)=\frac{b_{1}F(z)-a_{1}G(z)}{a_{0}b_{1}-b_{0}a_{1}}e^{\frac{\alpha z+\beta}{3}}=\frac{b_{1}-a_{1}\eta+\frac{b_{1}+a_{1}\eta}{\sqrt{3}}\wp'(h(z))}{2(a_{0}b_{1}-b_{0}a_{1})\wp(h(z))}e^{\frac{\alpha z+\beta}{3}}.
\end{equation}
By differentiating (\ref{eq-new8}) and noting that $(\wp')^{2}=4\wp^{3}-1$, $\wp''=6\wp^{2}$, we have
\begin{equation}\label{eq-new8-1}
\begin{aligned}
f'(z)=&\frac{\frac{2(b_{1}+a_{1}\eta)}{\sqrt{3}}\wp^{3}(h(z))h'(z)-(b_{1}-a_{1}\eta)\wp'(h(z))h'(z)+\frac{b_{1}+a_{1}\eta}{\sqrt{3}}h'(z)}{2(a_{0}b_{1}-b_{0}a_{1})\wp^{2}(h(z))}e^{\frac{\alpha z+\beta}{3}}\\
&+\frac{\alpha\left(b_{1}-a_{1}\eta+\frac{b_{1}+a_{1}\eta}{\sqrt{3}}\wp'(h(z))\right)}{6(a_{0}b_{1}-b_{0}a_{1})\wp(h(z))}e^{\frac{\alpha z+\beta}{3}}.
\end{aligned}
\end{equation}
Substituting (\ref{eq-new8}) and (\ref{eq-new8-1}) into the first equation of (\ref{eq-new1-1}), we have
\begin{equation}\label{eq-new8-2}
\begin{aligned}
&\Bigg((a_{0}b_{1}-b_{0}a_{1})\Big(1+\frac{\wp'(h(z))}{\sqrt{3}}\Big)\wp(h(z))-(a_{0}+\frac{a_{2}\alpha}{3})\Big(b_{1}-a_{1}\eta\\
&+\frac{b_{1}+a_{1}\eta}{\sqrt{3}}\wp'(h(z))\Big)\wp(h(z))
-a_{2}\Big(\frac{2(b_{1}+a_{1}\eta)}{\sqrt{3}}\wp^{3}(h(z))h'(z)\\
&-(b_{1}-a_{1}\eta)\wp'(h(z))h'(z)+\frac{b_{1}+a_{1}\eta}{\sqrt{3}}h'(z)\Big)\Bigg)\wp(h(z+c))\\
=&e^{\frac{\alpha c}{3}}a_{1}\left(b_{1}-a_{1}\eta+\frac{b_{1}+a_{1}\eta}{\sqrt{3}}\wp'(h(z+c))\right)\wp^{2}(h(z)).
 \end{aligned}
\end{equation}

Similar to Claim 2, we prove the following fact.

\textbf{Claim 3:} $\wp(h(a_{j,k}+c))=0$ only holds for at most finitely many $a_{j,k}$'s.

We assume, to the contrary, that there exists an infinite subsequence of  $\{a_{j,k}\}$, without loss of generality we may take $\{a_{j,k}\}$  itself,  such that $\wp(h(a_{j,k}+c))=0$. Thus we have
$(\wp')^{2}(h(a_{j,k}+c))=-1$.
 Differentiating (\ref{eq-new8-2}) and then setting $z=a_{j,k}$, we derive a relation
 \begin{equation}\label{new-8-3}
   a_2 h'(a_{j,k})h'(a_{j,k}+c)\left((a_{1}\eta-b_1)\wp'(h(a_{j,k}))+\frac{b_{1}+a_{1}\eta}{\sqrt{3}}\right)=0.
 \end{equation}
Note that $\wp'(h(a_{j,k}))=\pm i$, where $i$ is the imaginary unit. If $a_{2}\neq0$, the equation (\ref{new-8-3}) becomes
\begin{equation*}
  h'(a_{j,k})h'(a_{j,k}+c)\left(\pm(a_{1}\eta-b_1)i+\frac{b_{1}+a_{1}\eta}{\sqrt{3}}\right)=0.
\end{equation*}
Since $h$ is a nonconstant polynomial and $a_{j,k}$ is an infinite sequence, then there are infinitely many $a_{j,k}$ such that  $h'(a_{j,k})h'(a_{j,k}+c)\neq 0$. It follows that
\begin{equation}\label{two-cons-equ}
   \pm (a_{1}\eta-b_1)i+\frac{b_{1}+a_{1}\eta}{\sqrt{3}}=0.
\end{equation}
Let $t_j\ (j\geq 1)$ be the poles of $\wp(z)$ satisfying $|t_{j}|\rightarrow \infty$ as $j\rightarrow \infty$ and take $b_{j,k}\in\mathbb{C}$ satisfying $h(b_{j,k})=t_j$ for $k=1,...,\deg(h)$. Then there exists an integer $j_0$ such that when $j>j_0$, $b_{j,k}$ are simple zeros of $h(z)-t_j$  and $h(z+c)-t_j$ has only simple zeros. Thus, the unique term
$ 2(b_{1}+a_{1}\eta)\wp^{3}(h(z))h'(z)/{\sqrt{3}} $
with poles of multiplicity $6$ at $b_{j,k}\ (j>j_0)$ must vanish, that is, $b_{1}+a_{1}\eta=0$. Combining with (\ref{two-cons-equ}), we obtain $a_1=b_1=0$. This contradicts the assumption $a_0b_1-a_1b_0\not=0$.

If $a_{2}=0$, then we have $a_{1}b_{2}=0$ by the assumption $b_{1}a_{2}-a_{1}b_{2}=0$. Hence we distinguish two cases depending on whether $a_{1}$ is zero or not as follows.

If $a_{1}=0$, then we have $a_{0}b_{1}\neq0$ by the assumption $a_{0}b_{1}-a_{1}b_{0}\neq0$. Now, by using (\ref{eq-new1}), we can rewrite the equations (\ref{eq-new1-1}) as follows
\begin{equation}\label{eq-new9}
f(z)=\frac{1+\frac{\wp'(h(z))}{\sqrt{3}}}{2a_0\wp(h(z))}e^{\frac{\alpha z+\beta}{3}},\
b_{0}f(z)+b_{1}f(z+c)+b_{2}f'(z)=\frac{1-\frac{\wp'(h(z))}{\sqrt{3}}}{2\wp(h(z))}\eta e^{\frac{\alpha z+\beta}{3}}.
\end{equation}
Differentiating the first relation of (\ref{eq-new9}) and noticing that $(\wp')^{2}=4\wp^{3}-1$, $\wp''=6\wp^{2}$, it follows that
\begin{equation}\label{eq-new9-1}
\begin{aligned}
f'(z)=\frac{\frac{2\wp^{3}(h(z))h'(z)}{\sqrt{3}}-\wp'(h(z))h'(z)+\frac{h'(z)}{\sqrt{3}}}{2a_{0}\wp^{2}(h(z))}e^{\frac{\alpha z+\beta}{3}}
+\frac{\alpha(1+\frac{\wp'(h(z))}{\sqrt{3}})}{6a_{0}\wp(h(z))}e^{\frac{\alpha z+\beta}{3}}.
\end{aligned}
\end{equation}
Substituting $f$ and $f'$ into the second equation of (\ref{eq-new9}), we have
\begin{equation}\label{eq-new9-2}
\begin{aligned}
&\Bigg (a_{0}\eta\left(1-\frac{\wp'(h(z))}{\sqrt{3}}\right)\wp(h(z))-\left(b_{0}+\frac{b_{2}\alpha}{3}\right)\left(1+\frac{\wp'(h(z))}{\sqrt{3}}\right)\wp(h(z))\\
&-b_{2}\left(\frac{2\wp^{3}(h(z))h'(z)}{\sqrt{3}}-\wp'(h(z))h'(z)+\frac{h'(z)}{\sqrt{3}}\right)\Bigg)\wp(h(z+c))\\
 =&e^{\frac{\alpha c}{3}}b_{1}\left(1+\frac{\wp'(h(z+c))}{\sqrt{3}}\right)\wp^{2}(h(z)).
 \end{aligned}
\end{equation}
 Differentiating (\ref{eq-new9-2}) and then setting $z=a_{j,k}$, we derive a relation
 \begin{equation}\label{new-9-3}
   b_2 h'(a_{j,k})h'(a_{j,k}+c)\left(\wp'(h(a_{j,k}))-\frac{1}{\sqrt{3}}\right)=0.
 \end{equation}
Note that $\wp'(h(a_{j,k}))=\pm i$, where $i$ is the imaginary unit. If $b_{2}\neq0$, (\ref{new-9-3}) derives a contradiction.
Hence we have $b_2=0$. Now we can rewrite the second equation of (\ref{eq-new9}) as follows
\begin{equation*}
b_{0}f(z)+b_{1}f(z+c)=\frac{1-\frac{\wp'(h(z))}{\sqrt{3}}}{2\wp(h(z))}\eta e^{\frac{\alpha z+\beta}{3}}.
\end{equation*}
Substituting the first relation of (\ref{eq-new9}) into the above equation, we have
\begin{equation}\label{eq-new-9-6}
b_{1}e^{\frac{\alpha c}{3}}\Big(1+\frac{\wp'(h(z+c))}{\sqrt{3}}\Big)\wp(h(z))
 =\Big(a_{0}\eta-b_{0}-\frac{a_{0}\eta+b_{0}}{\sqrt{3}}\wp'(h(z))\Big)\wp(h(z+c)).
 \end{equation}
 Differentiating the above equation and then setting $z=a_{j,k}$, we derive a relation
 \begin{equation*}
 \begin{aligned}
  &b_{1}e^{\frac{\alpha c}{3}}\Big(1+\frac{\wp'(h(a_{j,k}+c))}{\sqrt{3}}\Big)\wp'(h(a_{j,k}))h'(a_{j,k})\\
 = & \Big(a_{0}\eta-b_{0}-\frac{a_{0}\eta+b_{0}}{\sqrt{3}}\wp'(h(a_{j,k}))\Big)\wp'(h(a_{j,k}+c))h'(a_{j,k}+c).
  \end{aligned}
 \end{equation*}
 Note that $\wp'(h(a_{j,k}))=\pm i$ and $\wp'(h(a_{j,k}+c))=\pm i$, where $i$ is the imaginary unit. The above equation immediately implies  one and only one of the following four situations
 \begin{equation}\label{eq-neww01}
 \begin{aligned}
  A_{1}h'(a_{j,k})&=B_{1}h'(a_{j,k}+c),\\
  A'_{1}h'(a_{j,k})&=B_{2}h'(a_{j,k}+c),\\
  A_{2}h'(a_{j,k})&=B_{1}h'(a_{j,k}+c),\\
  A'_{2}h'(a_{j,k})&=B_{2}h'(a_{j,k}+c),\\
  \end{aligned}
 \end{equation}
 where
 \begin{equation*}
   A_{1}=-A'_{1}= b_{1}e^{\frac{\alpha c}{3}}\left(1+\frac{i}{\sqrt{3}}\right), ~~B_{1}=a_{0}\eta-b_{0}-\frac{a_{0}\eta+b_{0}}{\sqrt{3}}i,
 \end{equation*}

 \begin{equation*}
  A_{2}=-A'_{2}= -b_{1}e^{\frac{\alpha c}{3}}\left(1-\frac{i}{\sqrt{3}}\right), ~~B_{2}=a_{0}\eta-b_{0}+\frac{a_{0}\eta+b_{0}}{\sqrt{3}}i.
 \end{equation*}
Obviously, $A_1, A_2$ are nonzero constants by the assumption  $a_{0}b_{1}\neq0$. Since $h$ is a nonconstant polynomial and $a_{j,k}$ is an infinite sequence, the relations (\ref{eq-neww01}) immediately yield functional equations
\begin{equation}\label{eq-neww02}
 \begin{aligned}
 A_{1}h'(z)&=B_{1}h'(z+c),\\
  A'_{1}h'(z)&=B_{2}h'(z+c),\\
  A_{2}h'(z)&=B_{1}h'(z+c),\\
  A'_{2}h'(z)&=B_{2}h'(z+c),\\
  \end{aligned}
 \end{equation}
which further mean that one of the following four equalities
\begin{equation*}
    A_1=B_1,\ -A_1=B_2,\ A_2=B_1,\ -A_2=B_2
\end{equation*}
holds by comparing the leading coefficient.

Thus we obtain $h'(z)=h'(z+c)$, which implies $h(z)=az+b$, where $a(\neq 0)$ and $b$ are constants. Now the equations
   $\wp (h(a_{j,k}))=0,\ \wp (h(a_{j,k}+c))=0$
become $ \wp (aa_{j,k}+b)=0,\ \wp (aa_{j,k}+b+ac)=0,$
 that is, $\{aa_{j,k}+b,aa_{j,k}+b+ac\}\subset\{z_j\}_{j=1}^\infty$. Hence we have
 \begin{equation*}
    (aa_{j,k}+b+ac)-(aa_{j,k}+b)=ac\in P_a\cup\{m\omega_1+n\omega_2:m,n\in\mathbb{Z}\},
 \end{equation*}
 since $\wp $ has only two distinct zeros in the parallelogram $P_{a}$ and $\wp$ is a function of double periods $\omega_1$ and $\omega_2$.

 If $ac=m\omega_1+n\omega_2$ for some $m,n\in\mathbb{Z}$, the equation (\ref{eq-new-9-6}) becomes
\begin{equation}\label{eq-neww05}
 \begin{aligned}
  \frac{a_{0}\eta+b_{0}+b_{1}e^{\frac{\alpha c}{3}}}{\sqrt{3}}\wp'(az+b)
  =a_{0}\eta-b_{0}-b_{1}e^{\frac{\alpha c}{3}}
  \end{aligned}
 \end{equation}
 since $\wp(az+b)=\wp(az+b+ac)$, and hence $a_{0}\eta+b_{0}+b_{1}e^{\frac{\alpha c}{3}}=0,\ a_{0}\eta-b_{0}-b_{1}e^{\frac{\alpha c}{3}}=0,$
because $\wp'(az+b)$ is a transcendental meromorphic function. We can obtain $a_{0}=0$, which is a contradiction.

When $ac\in P_a$, we rewrite (\ref{eq-new-9-6}) into the form
\begin{equation}\label{eq-neww06}
  \frac{a_{0}\eta-b_{0}-\frac{a_{0}\eta+b_{0}}{\sqrt{3}}\wp'(az+b)}{\wp(az+b)}
  =\frac{b_{1}e^{\frac{\alpha c}{3}}\Bigg(1+\frac{\wp'(az+b+ac)}{\sqrt{3}}\Bigg)}{\wp(az+b+ac)}.
\end{equation}
It is obvious that the function on the left-hand side of (\ref{eq-neww06}) has pole at $z=-\frac{b}{a}$, but the function on the right-hand side of (\ref{eq-neww06})  take a finite value. This leads to a contradiction.

Therefore, we must have $a_1\not=0$. By our assumptions $a_{2}=0$ and $b_{1}a_{2}-b_{2}a_{1}=0$, then $b_{2}=0$. Now, by using (\ref{eq-new1}), we can rewrite the equations (\ref{eq-new1-1}) as follows
\begin{equation}\label{eq-new11}
\begin{aligned}
  & a_{0}f(z)+a_{1}f(z+c)=\frac{1}{2}\frac{1+\frac{\wp'(h(z))}{\sqrt{3}}}{\wp(h(z))}e^{\frac{\alpha z+\beta}{3}},\\
  & b_{0}f(z)+b_{1}f(z+c)=\frac{\eta}{2}\frac{1-\frac{\wp'(h(z))}{\sqrt{3}}}{\wp(h(z))}e^{\frac{\alpha z+\beta}{3}}.
  \end{aligned}
\end{equation}
Solving the equations (\ref{eq-new11}), it follows that
\begin{equation}\label{eq-neww1}
 f(z)=\frac{b_{1}-a_{1}\eta+\frac{(b_{1}+a_{1}\eta)\wp'(h(z))}{\sqrt{3}}}{2(a_{0}b_{1}-a_{1}b_{0})\wp(h(z))}e^{\frac{\alpha z+\beta}{3}},
  f(z+c)=\frac{b_{0}-a_{0}\eta+\frac{(b_{0}+a_{0}\eta)\wp'(h(z))}{\sqrt{3}}}{2(a_{1}b_{0}-a_{0}b_{1})\wp(h(z))}e^{\frac{\alpha z+\beta}{3}},
\end{equation}
which derives an equation
\begin{equation}\label{eq-neww2}
\begin{aligned}
 &\left(b_{0}-a_{0}\eta+\frac{b_{0}+a_{0}\eta}{\sqrt{3}}\wp'(h(z))\right)\wp(h(z+c))\\
&= -e^{\frac{\alpha c}{3}}\left(b_{1}-a_{1}\eta +\frac{b_{1}+a_{1}\eta}{\sqrt{3}}\wp'(h(z+c))\right)\wp(h(z)).
\end{aligned}
\end{equation}
Differentiating the above equation and then setting $z=a_{j,k}$, we derive a relation
 \begin{equation*}
\begin{aligned}
 &\left(b_{0}-a_{0}\eta+\frac{b_{0}+a_{0}\eta}{\sqrt{3}}\wp'(h(a_{j,k}))\right)\wp'(h(a_{j,k}+c))h'(a_{j,k}+c)\\
= &-e^{\frac{\alpha c}{3}}\left(b_{1}-a_{1}\eta +\frac{b_{1}+a_{1}\eta}{\sqrt{3}}\wp'(h(a_{j,k}+c))\right)\wp'(h(a_{j,k}))h'(a_{j,k}).
\end{aligned}
\end{equation*}
Note that $\wp'(h(a_{j,k}))=\pm i$ and $\wp'(h(a_{j,k}+c))=\pm i$, where $i$ is the imaginary unit. The above equation immediately implies  one and only one of (\ref{eq-neww01})
 \begin{equation*}
   A_{1}=A'_{1}= b_{0}-a_{0}\eta+\frac{b_{0}+a_{0}\eta}{\sqrt{3}}i, ~~B_{1}=-e^{\frac{\alpha c}{3}}\left(b_{1}-a_{1}\eta +\frac{b_{1}+a_{1}\eta}{\sqrt{3}}i\right),
 \end{equation*}
 \begin{equation*}
  A_{2}=A'_{2}=-(b_{0}-a_{0}\eta-\frac{b_{0}+a_{0}\eta}{\sqrt{3}}i), ~~B_{2}=e^{\frac{\alpha c}{3}}\left(b_{1}-a_{1}\eta -\frac{b_{1}+a_{1}\eta}{\sqrt{3}}i\right).
 \end{equation*}
Note that $h$ is a nonconstant polynomial and $a_{j,k}$ is an infinite sequence, immediately yields the corresponding functional equation (\ref{eq-neww02}),
which further mean that one of the following four equalities
\begin{equation}\label{eq-A1}
     A_1=B_1,\ A_1=B_2,\ A_2=B_1,\ A_2=B_2
\end{equation}
holds by comparing the leading coefficient.

In the following, we only consider the case $A_1=B_1$ in (\ref{eq-A1}), and omit the others due to the similarity
of their arguments.
Note that $\wp'(h(a_{j,k}))=\pm i$. W. l. o. g., we assume $\wp'(h(a_{j,k}))= i$. Now, if $A_1$ is zero. Further, $a_{j,k}$ is also a zero of
\begin{equation*}
  b_{0}-a_{0}\eta+\frac{b_{0}+a_{0}\eta}{\sqrt{3}}\wp'(h(z)).
\end{equation*}
Assume that $a_{j,k}$ is a zero of $\wp(h(z))$ with multiplicity $l$. It follows from (\ref{eq-neww1}) that $a_{j,k}$ may be a pole of $f(z)$ and $f(z+c)$ of multiplicity $\leq l-1$. However, it follows from (\ref{eq-new11}) that $a_{j,k}$ is a pole of
$a_{0}f(z)+a_{1}f(z+c)$
with multiplicity $l$. This is a contradiction.

Hence we have $A_1\not=0$. Similarly, if $\wp'(h(a_{j,k}))=- i$, we also obtain $A_2\not=0$. Any way, we also have $B_1\not=0, B_2\not=0$.
Any way, we also have $B_1\not=0, B_2\not=0$. Thus we obtain $h'(z)=h'(z+c)$ again from (\ref{eq-neww02}), that is, $h(z)=az+b$, where $a(\neq 0)$ and $ b$ are constants.
Next we can derives a contradiction according to the arguments between (\ref{eq-neww02}) and (\ref{eq-neww06}).
Hence Claim 3 is proved.

Based on Claim 3, we can obtain a contradiction according to the proof of Case 1 after Claim 2.

Hence $F=c_0$ and $G=c_1$ are constants such that $c_0^3+c_1^3=1$. By solving (\ref{eq-new8}), we obtain a solution
\begin{equation*}
  f(z)=\frac{b_{1}c_{0}-a_{1}c_{1}}{a_{0}b_{1}-b_{0}a_{1}}e^{\frac{\alpha z+\beta}{3}}.
\end{equation*}
This is just Case $2$ in Theorem  \ref{mai-1.4}.

\textbf{Case 3:} $b_{1}a_{2}-a_{1}b_{2}\neq0$.

 From (\ref{eq-new1-1}) and (\ref{eq-new1}), we have
\begin{equation*}
  (b_{1}a_{2}-a_{1}b_{2})f'(z)+(a_{0}b_{1}-a_{1}b_{0})f(z)=(b_{1}F(z)-a_{1}G(z))e^{\frac{\alpha z+\beta}{3}}.
\end{equation*}
By using the theory of first-order linear differential equations, the above equation has the solution
\begin{equation}\label{eq-new12}
  f(z)=e^{-\frac{a_{0}b_{1}-a_{1}b_{0}}{b_{1}a_{2}-a_{1}b_{2}}z}\left(\int\frac{b_{1}-a_{1}\eta+\frac{b_{1}+a_{1}\eta}{\sqrt{3}}\wp'(h(z))}{2(b_{1}a_{2}-a_{1}b_{2})\wp(h(z))}
  e^{(\frac{a_{0}b_{1}-a_{1}b_{0}}{b_{1}a_{2}-a_{1}b_{2}}+\frac{\alpha}{3})z+\frac{\beta}{3}}dz+C\right),
\end{equation}
where $C$ is a constant.

If $b_{1}+a_{1}\eta\neq 0$, letting $t_j\ (j\geq 1)$ be the poles of $\wp(z)$, then we have
$\wp(z)=g_j(z)(z-t_{j})^{-2},$
where $g_j$ is a holomorphic function in a neighborhood of $t_j$ with $g_j(t_{j})\neq0$, and hence
\begin{equation*}
\begin{aligned}
\frac{b_{1}-a_{1}\eta+\frac{b_{1}+a_{1}\eta}{\sqrt{3}}\wp'(z)}{2(b_{1}a_{2}-a_{1}b_{2})\wp(z)}=&-\frac{b_{1}+a_{1}\eta}{\sqrt{3}(b_{1}a_{2}-a_{1}b_{2})}\frac{1}{z-t_{j}}
+\frac{(b_{1}-a_{1}\eta)(z-t_{j})^{2}}{2(b_{1}a_{2}-a_{1}b_{2})g_j(z)}\\
&+\frac{b_{1}+a_{1}\eta}{2\sqrt{3}(b_{1}a_{2}-a_{1}b_{2})}\frac{g_j'(z)}{g_j(z)}\\
=&-\frac{b_{1}+a_{1}\eta}{\sqrt{3}(b_{1}a_{2}-a_{1}b_{2})}\frac{1}{z-t_{j}}+O(1)
\end{aligned}
\end{equation*}
as $z\rightarrow t_{j}$. Setting $z=h(z)$, we derive a relation
\begin{equation}\label{equ-1}
\begin{aligned}
\frac{b_{1}-a_{1}\eta+\frac{b_{1}+a_{1}\eta}{\sqrt{3}}\wp'(h(z))}{2(b_{1}a_{2}-a_{1}b_{2})\wp(h(z))}=-\frac{b_{1}+a_{1}\eta}{\sqrt{3}(b_{1}a_{2}-a_{1}b_{2})}\frac{1}{h(z)-t_{j}}+O(1)
\end{aligned}
\end{equation}
as $h(z)\to t_j$. Note that the equation $h'(z)=0$ has only finitely many solutions, and $\bigcup \limits^{\infty}_{j=1}h^{-1}(t_{j})$ is an infinite set. There exist an integer $j$ and a point $z'\in h^{-1}(t_{j})$ such that $h'(z')\neq 0$, that is, $z'$ is a simple zero of $h(z)-t_j$. It follows from (\ref{equ-1}) that $f(z)$ has a logarithmic singular point $z'$, which contradicts to the assumption  that $f(z)$ is a meromorphic function.

If $b_{1}+a_{1}\eta=0$, we have $b_{1}-a_{1}\eta\neq 0$. Otherwise, if $b_{1}-a_{1}\eta=0$, that is, $b_1=0$ and hence $a_1=0$, which is contracted to the assumption $b_1a_2-a_1b_2\not=0$. Now, we can rewrite the equation (\ref{eq-new12}) as follows
\begin{equation*}
  f(z)=e^{-\frac{a_{0}b_{1}-a_{1}b_{0}}{b_{1}a_{2}-a_{1}b_{2}}z}\left(\int\frac{b_{1}-a_{1}\eta}{2(b_{1}a_{2}-a_{1}b_{2})\wp(h(z))}
  e^{(\frac{a_{0}b_{1}-a_{1}b_{0}}{b_{1}a_{2}-a_{1}b_{2}}+\frac{\alpha}{3})z+\frac{\beta}{3}}dz+C\right).
\end{equation*}
Similarly, there exist an integer $j$ and a point $z''\in h^{-1}(z_{j})$ such that $h'(z'')\neq 0$,
that is,  $f(z)$ has a logarithmic singular point $z''$ since $\wp $ has only simple zeros. We obtain a contradiction again.

Hence $F=c_0$ and $G=c_1$ are constants such that $c_0^3+c_1^3=1$. Thus (\ref{eq-new12}) admits a solution
\begin{equation*}
  f(z)=\frac{3(b_{1}c_{0}-a_{1}c_{1})}{(b_{1}a_{2}-a_{1}b_{2})(\alpha-3D)}e^{\frac{\alpha z+\beta}{3}}+Ce^{Dz},
\end{equation*}
where $D=\frac{a_{1}b_{0}-a_{0}b_{1}}{b_{1}a_{2}-a_{1}b_{2}}$ is a constant. This is just Case $3$ in Theorem  \ref{mai-1.4}. The proof of Theorem \ref{mai-1.4} is completed.




\begin{thebibliography}{99}
\addcontentsline{toc}{chapter}{Bibliography}
\bibitem{15} I.N. Baker,  On a class of meromorphic functions. Pro. Amer. Math. Soc. 17 (4) (1966) 819--822.

\bibitem{22} S.B. Bank, J.K. Langley, {\em On the value distribution theory of elliptic functions.} Monatsh. Math. 98 (1) (1984) 1--20.

\bibitem{24} W. Bergweiler, {\em Order and lower order of composite meromorphic functions}. Michigan Math. J. 36 (1989) 135--146.

\bibitem{10} D.C. Chang, B.Q. Li, {\em Description of entire solutions of eiconal type equations}. Canad. Math. Bull. 55 (2012) 249--259.

\bibitem{25} Y.M. Chiang, S. Feng, {\em On the Nevanlinna characteristic of $f(z +\eta)$ and difference equations in the
complex plane}. Ramanujan J. 16 (2008) 105--129.

\bibitem{1} A. Edrei, W.H.J. Fuchs, {\em On the zeros of $f(g(z))$ where $f$ and $g$ are entire functions}. J. Anal. Math. 12 (1964) 243--255.

\bibitem{11} F. Gross,  {\em On the equation $ f^{n}+ g^{n}=1$.} I. Bull. Amer. Math. Soc. 72 (1) (1966) 86--88.

\bibitem{12} F. Gross, {\em On the functional equation $f^{n}+ g^{n}= h^{n}$.} Amer. Math. Monthly 73 (10) (1966) 1093--1096.

\bibitem{13} F. Gross, {\em On the equation $ f^{n}+ g^{n}=1$.} II. Bull. Amer. Math. Soc. 74 (1968) 647--767.

\bibitem{35} G.G Gundersen, {\em Research questions on meromorphic functions and complex differential equations}. Comput. Methods Funct. Theory 17 (2017) 195--209.

\bibitem{40} G.G. Gundersen, {\em Meromorphic solutions of $f^{5}+g^{5}+h^{5}\equiv1$.} Complex Var. Theory Appl. 43 (2001) 293--298.

\bibitem{41} G.G. Gundersen, {\em Meromorphic solutions of $f^{6}+g^{6}+h^{6}\equiv1$.} Analysis (Munich). 18 (1998) 285--290.


\bibitem{19} Q. Han, F. L\"u, {\em On the functional equation $f^{n} (z)+ g^{n} (z)= e^{\alpha z+\beta}$. } arXiv:1612.06842v1 (2016).

\bibitem{14} Q. Han, F. L\"u, , {\em On the {F}ermat-type equation {$f^3(z)+f^3(z+c)=1$}.} Aequationes math. 91 (1) (2017) 129--136.

\bibitem{20} Q. Han, {\em On complex analytic solutions of the partial differential equation $(u_{z_{1}} )^{m} + (u_{z_{2}} )^{n} = u^{m}$.} Houston J. Math. 35 (2009) 277--289.

\bibitem{32} R.G. Halburd, R.J. Korhonen, {\em Nevanlinna theory for the difference operator}. Ann. Acad. Sci. Fenn. Math. 31 (2006) 463--478.

\bibitem{31} R.G. Halburd, R.J. Korhonen, {\em Difference analogue of the lemma on the logarithmic derivative with applications to difference equations.} Journal of Mathematical Analysis and Applications 314 (2) (2006) 477-487.

\bibitem{43} R.G. Halburd, R.J. Korhonen, {\em Meromorphic solutions of difference equations, integrability
           and the discrete P ainlev$\acute{e}$ equations.} J. Phys. A. 40 (2007) 1--38.

\bibitem{28} W. Hayman, {\em Waring's Problem f\"ur analytische Funktionen.} Bayer. Akad. Wiss. Math.-Natur. Kl. Sitzungsber. 1984 (1985) 1--13.

\bibitem{26} P.C. Hu, P. Li, C.C. Yang, {\em Unicity of meromorphic mappings.} Adv. Complex Anal. Appl. Kluwer Academic Publishers, 2003.

\bibitem{27} K. Ishizaki, {\em A note on the functional equation $f^{n}+g^{n}+ h^{n}=1$ and some complex differential equations.} Comput. Methods Funct. Theory 2 (1) (2003) 67--85.


\bibitem{45} I. Laine, {\em Nevanlinna Theory and Complex Differential Equations.} Walter de Gruyter, Berlin-New York (1993).

\bibitem{16} K. Liu, L.Z. Yang, X.L. Liu, {\em Existence of entire solutions of nonlinear difference equations.} Czech. Math. J. 61 (2011) 565--576.

\bibitem{17} K. Liu, T.B. Cao, H.Z. Cao, {\em Entire solutions of Fermat type differential-difference equations.} Arch. Math. (Basel) 99 (2012) 147--155.

\bibitem{2} B.Q. Li, {\em On meromorphic solutions of $f^{2} + g^{2} = 1$.} Math. Z. 258 (2008) 763--771.

\bibitem{3} B.Q. Li, Z. Ye, {\em On meromorphic solutions of $f^{3} + g^{3} = 1$.} Arch. Math. (Basel) 90 (2008) 39--43.

\bibitem{4} B.Q. Li, {\em Entire solutions of $(u_{z_{1}} )^{m} + (u_{z_{2}} )^{n} = e^{g}$.} Nagoya Math. J. 178 (2005) 151-162.

\bibitem{5} B.Q. Li, {\em Entire solutions of eiconal type equations.} Arch. Math. (Basel) 89 (2007) 350--357.

\bibitem{18} P. Montel, {\em Le\c{c}ons sur les familles normales de fonctions analytiques et leurs applications.} Collect. Borel, Paris, 1927, pp. 135--136.

\bibitem{21} R. Nevanlinna, {\em Analytic Functions.} Springer-Verlag, New York-Berlin, 1970.

\bibitem{29} T.W. Ng, S.K. Yeung, {\em Entire holomorphic curves on a Fermat surface of low degree}. arXiv:1612.01290 (2016).

\bibitem{23} S. Shimomura, {\em Entire solutions of a polynomial difference equation.} J. Fac. Sci. Univ. Tokyo Sect. IA Math. 28 (1981) 253--266.

\bibitem{7} C.C. Yang, {\em A generalization of a theorem of P. Montel on entire functions}. Proc. Amer. Math. Soc. 26 (1970) 332--334.

\bibitem{8} C.C. Yang, I. Laine, {\em On analogies between nonlinear difference and differential equations.} Proc. Japan Acad., Ser. A. 86 (2010) 10--14.

\bibitem{9} C.C. Yang, P. Li, {\em On the transcendental solutions of a certain type of nonlinear differential equations.} Arch. Math. 82 (2004) 442--448.

\bibitem{42} L.Z Yang, Z. Zhong, J.L Zhang, {\em Non-existence of meromorphic solutions of a Fermat type functional equation.} Aequationes math. 76 (2008) 140--150.

\end{thebibliography}
\end{document}